\documentclass[11pt]{amsart}

\usepackage{amsmath}
\usepackage{amssymb}
\usepackage{amscd}
\usepackage{epsfig}

 \newcommand{\N}{\mathbb N}
 \newcommand{\bee}{\begin{equation}}
 \newcommand{\eee}{\end{equation}}
  
 \newcommand{\Lb}{\mbox {\boldmath ${\Lambda}$}}

 \newcommand{\Lbs}{\mbox{\scriptsize\boldmath ${\Lambda}$}}
 
\def\Nk{{\mathcal N}}

 \newcommand{\Lamt}{\mbox{\tiny ${\Lambda}$}}
 \newcommand{\Pb}{\mbox {\bf P}}
 \newcommand{\Pbs}{\mbox {\scriptsize{\bf P}}}

\newcommand{\be}{\begin{eqnarray}}
\newcommand{\ee}{\end{eqnarray}}
\newcommand{\supp}{\mbox{\rm supp}}

\newcommand{\freq}{\mbox{\rm freq}}
\newcommand{\Vol}{\mbox{\rm Vol}}

\newcommand{\eps}{{\mbox{$\epsilon$}}}
\newcommand{\e}{{\varepsilon}}

\newcommand{\gam}{{\gamma}}

\newcommand{\R}{{\mathbb R}}
\newcommand{\Q}{{\mathbb Q}}
\newcommand{\Z}{{\mathbb Z}}
\newcommand{\C}{{\mathbb C}}
\newcommand{\Nat}{{\mathbb N}}

\newcommand{\Ak}{{\mathcal A}}
\newcommand{\Hk}{{\mathcal H}}

\newcommand{\Dk}{{\mathcal D}}

\newcommand{\Pk}{{\mathcal P}}

\newcommand{\Cant}{{\mathcal C}}

\newcommand{\Gk}{{\mathcal G}}

\newcommand{\Sk}{{\mathcal S}}
\newcommand{\Tk}{{\mathcal T}}

\newcommand{\dist}{\mbox{\rm dist}}

\newcommand{\Lam}{{\Lambda}}

\newcommand{\lam}{\lambda}

\newcommand{\om}{\omega}

 \newtheorem{theorem}{Theorem}[section]
 \newtheorem{lemma}[theorem]{Lemma}
 \newtheorem{prop}[theorem]{Proposition}
 \newtheorem{cor}[theorem]{Corollary}

 \newtheorem{defi}[theorem]{Definition}
 
 \newtheorem{remark}[theorem]{Remark}

\numberwithin{equation}{section}

\begin{document}

\title[Meyer property for substitution Delone sets]{Pure point diffractive substitution Delone sets have the Meyer property}
 
\author{Jeong-Yup  Lee}
\address{Jeong-Yup  Lee, Department of Mathematics and Statistics,
Box 3045 STN CSC, University of Victoria, Victoria, BC, V8W 3P4, Canada}
\email{jylee@math.uvic.ca}

\author{Boris Solomyak}
\address{Boris Solomyak, Box 354350, Department of Mathematics,
University of Washington, Seattle WA 98195}
\email{solomyak@math.washington.edu}

\date{\today}

 \thanks{2000 {\em Mathematics Subject Classification:} Primary 52C23, 
Secondary 37B50 
\\ \indent
{\em Key words and phrases:} quasicrystal, Delone set, Meyer set,
substitution, tiling dynamical system \\ \indent 
The first author
acknowledges support from
the NSERC post-doctoral fellowship and thanks to
the University of Washington and the University of Victoria for being the host
universities of the fellowship.
\\ \indent
The second author
is grateful to the Weizmann Institute of Science where he was a Rosi and
Max Varon Visiting Professor when this work was completed. He was also
supported in part by NSF grant DMS 0355187.}

\begin{abstract}
We prove that a primitive substitution Delone set, which is pure point
diffractive, is a Meyer set. This answers a question of J. C. Lagarias.
We also show that for primitive substitution Delone sets, being a
Meyer set is equivalent to having a relatively dense set of Bragg peaks.
The proof is based on tiling dynamical systems and the connection
between the diffraction and dynamical spectra.
\end{abstract}

\maketitle

\thispagestyle{empty}

\section{Introduction}

The discovery of quasicrystals in the 1980's inspired a lot of research
in the area of ``aperiodic order'' and ``mathematical quasicrystals.''
Roughly speaking, physical quasicrystals are aperiodic structures which
exhibit sharp bright
spots (called Bragg peaks) in their $X$-ray diffraction pattern.
The presence of Bragg peaks indicates the presence of ``long-range order"
in the structure.
A mathematical idealization of a large set of atoms is a discrete set
in $\R^d$. The most general class of sets modeling solids 
is the class of {\em Delone sets}, that is, subsets of
$\R^d$ which are relatively dense and uniformly discrete. Usually
some additional assumptions are made. A Delone set $\Lam$ is
of {\em finite local
complexity} (or ``finite type'') if $\Lam-\Lam$ is closed
and discrete, which is equivalent to having finitely many local patterns,
up to translations, see \cite{Lag99}. 
Another common assumption is {\em repetitivity}, which means that every
pattern of a Delone set (and not just individual points) occurs relatively
densely in space.
This is still not enough for
long-range order, since a repetitive Delone set of finite local complexity
may fail to have any
Bragg peaks. The Delone set $\Lam$ is said to be a {\em Meyer set} if
$\Lam-\Lam$ is uniformly discrete. 
Meyer sets were introduced (under the name of ``harmonious sets'') 
in 1969-1970 by Y. Meyer \cite{Meyer}
in the context of harmonic analysis. In the last ten years their 
importance in the theory of long-range aperiodic order has been revealed
in many investigations, see e.g.\ \cite{RVM97, LMS2, Lee, BLM}. 

The mathematical concept of diffraction spectrum is based on the
Fourier transform of the autocorrelation measure, see \cite{Hof1,Hof2}. Under
certain conditions, this Fourier transform is a measure (called
{\em diffraction measure}) on $\R^d$, whose
discrete component corresponds to the Bragg peaks. A Delone set $\Lam$
is said to be {\em pure point diffractive} (or ``perfectly diffractive,''
or a ``Patterson set'' \cite{Lag00}) if the diffraction measure is pure
point (pure discrete).  There is another
notion of spectrum, which comes from Ergodic Theory via a dynamical
system associated with the Delone set. As shown by Dworkin \cite{Dw}
(see also \cite{LMS1,Gouere,BL}), there is a close connection between the
two notions of spectra.
 
In his survey on mathematical quasicrystals,
J. C. Lagarias raised the following problem \cite[Problem 4.10]{Lag00}.
{\em Let $\Lam$ be a Delone set of finite type which is repetitive. If $\Lam$ is pure point diffractive,
must it be a Meyer set?} We do not have an answer for this question, but we
solve the following special case.

\medskip

{\bf \cite[Problem 4.11]{Lag00}.} {\em Suppose that $\Lam$ is a primitive
self-replicating Delone set of finite type. If $\Lam$ is pure point diffractive,
must $\Lam$ be a Meyer set?}

\medskip
At this point, we
just mention that a primitive self-replicating Delone set, roughly speaking,
corresponds to the set of ``control points'' of a self-affine tiling.
In this paper we refer to 
it as a representable primitive substitution Delone set.
Precise definitions on representable primitive substitution Delone sets are given in the next section.

Our main result, Th.\,\ref{Ppd-onPointset-Meyer} answers this question affirmatively. This result is applicable to \cite{LMS2, Lee} in which the Meyer condition is additionally assumed to understand the structure of pure point diffractive point sets.

In fact, the condition of being pure point diffractive may be weakened.
We only need the fact that the set of Bragg peaks is relatively dense
in the entire space (this holds in the case of pure point diffractive set).
This condition turns out to be necessary and sufficient for the Meyer
property on the class of substitution Delone sets (see Th.\,\ref{Bragg-relativelyDense-Meyer}). 

The proof of the implication (for substitution Delone sets)
$$
\mbox{ relatively dense set of Bragg peaks}\ \Rightarrow\ \mbox{ Meyer set}
$$
relies on the theory of tiling dynamical systems developed in 
\cite{soltil} and the connection between substitution Delone sets,
substitution Delone set families, and self-affine tilings, studied in
\cite{lawa,LMS2}. The second key ingredient is a generalization of 
classical results by Pisot in Diophantine approximation, due to 
K\"ornyei \cite{Korn} and Mauduit \cite{Maud}. 
The relevance of PV-numbers (Pisot-Vijayaraghavan numbers) for the
Meyer set property was already pointed out by Meyer \cite{Meyer}.
We show that the expanding linear map associated with our
substitution Delone set satisfies the ``Pisot family'' condition
(this is 
essentially proved in \cite{Rob} based on \cite{soltil}), and we obtain
some extra information about the set of translation vectors between
tiles of the same type. The last ingredient is a generalization of the
well-known ``Garsia Lemma'' \cite[Lemma 1.51]{gar} (obtained independently by other
authors as well), which implies that the set of polynomials of arbitrary degree
with integer coefficients bounded by a uniform constant, evaluated
at a PV-number, yields a uniformly discrete set.

Now we can state our main result.

\begin{theorem} \label{th-main}
If $\Lam$ is a representable primitive
substitution Delone set of finite local complexity (FLC) such that
the Bragg peaks are relatively dense in $\R^d$, then $\Lam$ is a Meyer set.
\end{theorem}

We note that the converse is also true by a theorem of Strungaru \cite{Nik}:
if $\Lam$ is a Meyer set, then 
the Bragg peaks are relatively dense.

\begin{cor} \label{cor1}
If $\Lam$ is a representable primitive
substitution Delone set of finite local complexity which is pure point
diffractive, then $\Lam$ is a Meyer set.
\end{cor}

This resolves \cite[Problem 4.11]{Lag00} (it follows from the
context of \cite{Lag00} that FLC is implicitly assumed).

\section{Preliminaries}

\subsection{Substitution Delone multisets and tilings}

\noindent

A {\em multiset \footnote{Caution: In \cite{lawa}, 
the word multiset refers to a set with multiplicities.}} 
or {\em $m$-multiset} in $\R^d$ is a 
subset $\Lb = \Lam_1 \times \dots \times \Lam_m 
\subset \R^d \times \dots \times \R^d$ \; ($m$ copies)
where $\Lam_i \subset \R^d$. We also write 
$\Lb = (\Lam_1, \dots, \Lam_m) = (\Lam_i)_{i\le m}$.
Recall that a Delone set is a relatively dense and uniformly discrete 
subset of $\R^d$.
We say that $\Lb=(\Lambda_i)_{i\le m}$ is a {\em Delone multiset} in $\R^d$ if
each $\Lambda_i$ is Delone and $\supp(\Lb):=\bigcup_{i=1}^m \Lambda_i 
\subset \R^d$ is Delone.

Although $\Lb$ is a product of sets, it is convenient to think
of it as a set with types or colors, $i$ being the
color of points in $\Lambda_i$. 
A {\em cluster} of $\Lb$ is, by definition,
a family $\Pb = (P_i)_{i\le m}$ where $P_i \subset \Lambda_i$ is 
finite for all $i\le m$. For a bounded set $A \subset \R^d$, let
$A \cap \Lb := (A \cap \Lam_i)_{i \le m}$. 
There is a natural translation $\R^d$-action on the set of Delone multisets and their clusters
in $\R^d$. The translate of a cluster $\Pb$ by $x \in \R^d$ is $x + \Pb = (x+P_i)_{i \le m}$.
We say that two clusters $\Pb$ and $\Pb'$ are translationally equivalent if $\Pb = x + \Pb'$, i.e. 
$P_i = x + P_i'$ for all $i \le m$, for some $x \in \R^d$. We write $B_R(y)$ for the closed ball of radius $R$ centered at $y$.

\begin{defi} \label{def-flc}
{\em A Delone multiset $\Lb$ has {\em finite local complexity
(FLC)} if for every $R>0$ there exists a finite set $Y\subset \supp(\Lb)=
\bigcup_{i=1}^m \Lam_i$  such that 
$$
\forall x\in \supp(\Lb),\ \exists\, y\in Y:\ 
B_R(x) \cap \Lb = (B_R(y) \cap \Lb) + (x-y).
$$
In plain language, for each radius $R > 0$ there are only finitely many 
translational classes of clusters whose support lies in some ball of 
radius $R$.}
\end{defi}

\begin{defi}
{\em A Delone set $\Lambda$ is called a {\em Meyer set} if $\Lambda - \Lambda$ is 
uniformly discrete.}
\end{defi}

For a cluster $\Pb$ and a bounded set $A\subset \R^d$ denote
$$
L_{\Pbs}(A) = \sharp\{x\in \R^d:\ x+\Pb \subset A\cap \Lb\},
$$
where $\sharp$ means the cardinality.
In plain language, $L_{\Pbs}(A)$ is the number of translates of $\Pb$
contained in $A$, which is clearly finite.
For a bounded set $F \subset \R^d$ and $r > 0$, let $(F)^{+r}:=
\{x \in \R^d:\,\dist(x,F) \le r\}$ denote the $r$-neighborhood of $F$.
A {\em van Hove sequence} for $\R^d$ is a sequence 
$\mathcal{F}=\{F_n\}_{n \ge 1}$ of bounded measurable subsets of 
$\R^d$ satisfying
\be \label{Hove}
\lim_{n\to\infty} \Vol((\partial F_n)^{+r})/\Vol(F_n) = 0,~
\mbox{for all}~ r>0.
\ee

\begin{defi} \label{def-ucf}
Let $\{F_n\}_{n \ge 1}$ be a van Hove sequence.
The Delone multiset $\Lb$ has {\em uniform cluster frequencies} (UCF)
(relative to $\{F_n\}_{n \ge 1}$) if for any non-empty cluster $\Pb$, the limit
$$
\freq(\Pb,\Lb) = \lim_{n\to \infty} \frac{L_{\Pbs}(x+F_n)}{\Vol(F_n)} \ge 0
$$
exists uniformly in $x\in \R^d$.
\end{defi}

A linear map $Q : \R^d \rightarrow \R^d$ is {\em expansive}
if its every eigenvalue lies outside the unit circle. 

\begin{defi} \label{def-subst-mul}
{\em $\Lb = (\Lam_i)_{i\le m}$ is called a {\em
substitution Delone multiset} if $\Lb$ is a Delone multiset and
there exist an expansive map
$Q:\, \R^d\to \R^d$ and finite sets $\Dk_{ij}$ for $i,j\le m$ such that
\be \label{eq-sub}
\Lambda_i = \bigcup_{j=1}^m (Q \Lambda_j + \Dk_{ij}),\ \ \ i \le m,
\ee
where the unions on the right-hand side are disjoint.}
\end{defi}

For any given substitution Delone multiset $\Lb = (\Lambda_i)_{i \le m}$, 
we define $\Phi_{ij} = \{ f : x \mapsto Qx + a \, : \,a \in \Dk_{ij}\}$.
Then $\Phi_{ij}(\Lam_j) = Q \Lam_j + \Dk_{ij}$, where $i \le m$. We define $\Phi$ an $m \times m$ array for which each entry is $\Phi_{ij}$, and call $\Phi$ a {\em matrix function system (MFS)} for the substitution. 
For any $k \in \Z_+$ and $x \in \Lam_j$ with $j \le m$, we let
$\Phi^k (x) =  \Phi^{k-1}((\Phi_{ij}(x))_{i \le m})$.
                                                                                                                                        
We say that the
substitution Delone multiset $\Lb$ is {\em primitive} if the corresponding
substitution matrix $S$, with $S_{ij}= \sharp (\Dk_{ij})$, is primitive, i.e. there is an $l > 0$ for which $S^l$ has no zero entries.

We say that a Delone set $\Lam$ is a {\em substitution Delone set} if
there is a substitution Delone multiset $\Lb = (\Lam_i)_{i\le m}$ such that
$\Lam = \bigcup_{i=1}^m \Lam_i$. The Delone set $\Lam$ is said to be primitive if
the substitution Delone multiset $\Lb$ can be chosen primitive.

\medskip

Next we briefly review the basic definitions of tilings and substitution tilings.
We begin with a set of types (or colors) $\{1,\ldots,m\}$, 
which we fix once and for all. 
A {\em tile} in $\R^d$ is defined as a pair $T=(A,i)$ where $A=\supp(T)$ 
(the support of $T$) is a compact
set in $\R^d$ which is the closure of its interior, and 
$i=l(T)\in \{1,\ldots,m\}$
is the type of $T$. We let $g+T = (g+A,i)$ for $g\in \R^d$. We say that
a set $P$ of tiles is a {\em patch} if the number of tiles in $P$ is 
finite and the tiles of $P$ have mutually disjoint
interiors (strictly speaking, we have to say ``supports of tiles,'' but this
abuse of language should not lead to confusion). 
A tiling of $\R^d$ is a set $\Tk$ of tiles such that 
$\R^d = \cup \{\supp(T) : T \in \Tk\}$ and distinct tiles have disjoint 
interiors.
Given a tiling $\Tk$, finite sets of tiles of $\Tk$ are called 
$\Tk$-patches.

We define FLC and UCF for tilings 
in the same way as the corresponding properties for Delone multisets.

We always assume that any two $\Tk$-tiles with the same color are translationally equivalent.
(Hence there are finitely many $\Tk$-tiles up to translation.)

\begin{defi}\label{def-subst}
{\em Let $\Ak = \{T_1,\ldots,T_m\}$ be a finite set of tiles in $\R^d$
such that $T_i=(A_i,i)$; we will call them {\em prototiles}.
Denote by $\Pk_{\Ak}$ the set of
patches made of tiles each of which is a translate of one of $T_i$'s.
We say that $\omega: \Ak \to \Pk_{\Ak}$ is a {\em tile-substitution} (or simply
{\em substitution}) with
expansive map $Q$ if there exist finite sets $\Dk_{ij}\subset \R^d$ for
$i,j \le m$, such that
\begin{equation}
\om(T_j)= 
\{u+T_i:\ u\in \Dk_{ij},\ i=1,\ldots,m\} \ \ \  \mbox{for} \  j\le m,
\label{subdiv}
\end{equation}
with
$$
Q A_j  = \bigcup_{i=1}^m (\Dk_{ij}+A_i).
$$
Here all sets in the right-hand side must have disjoint interiors;
it is possible for some of the $\Dk_{ij}$ to be empty.}
\end{defi}

The substitution (\ref{subdiv}) is extended to all translates of prototiles by
$\om(x+T_j)= Q x + \om(T_j)$, and to patches and tilings by
$\om(P)=\cup\{\om(T):\ T\in P\}$.
The substitution $\om$ can be iterated, producing larger and larger patches
$\om^k(T_j)$. To the substitution $\om$ we associate its $m \times m$ 
substitution matrix $S$, with $S_{ij}:=\sharp (\Dk_{ij})$.
The substitution $\om$ is called {\em primitive}
if the substitution matrix $S$ is primitive.
We say that $\Tk$ is a fixed point of a substitution if $\om(\Tk) = \Tk$.

For each primitive substitution Delone multiset $\Lb$ (\ref{eq-sub}) one can set up  
an {\em adjoint system} of equations
\be \label{eq-til}
Q A_j = \bigcup_{i=1}^m (\Dk_{ij} + A_i),\ \ \ j \le m.
\ee
From Hutchinson's Theory (or rather, its generalization to the 
``graph-directed''
setting), it follows that (\ref{eq-til}) always has a unique solution 
for which $\Ak = \{A_1, \dots, A_m\}$ is 
a family of non-empty compact sets of $\R^d$ 
(see for example \cite{BM1}, Prop.\,1.3).
It is proved in \cite[Th.\,2.4 and Th.\,5.5]{lawa} that if $\Lb$ is a primitive
substitution Delone multiset, then all the sets $A_i$ from (\ref{eq-til})
have non-empty interiors and, moreover, each $A_i$ is the closure of 
its interior.

\begin{defi} {\em A Delone multiset $\Lb = (\Lam_i)_{i \le m}$ is called
{\em representable} (by tiles) for a tiling if there exists a set of prototiles 
$\Ak = \{T_i : i\le m\} $
so that
\be\label{eq-1}
\Lb + \Ak := \{x + T_i :\ x\in \Lambda_i,\ i \le m\} \ \ \ \mbox{is a tiling of}\ \ 
\R^d,
\ee
that is, $\R^d = \bigcup_{i\le m} \bigcup_{x\in \Lambda_i} (x + A_i)$ where $T_i = (A_i,i)$ for 
$i \le m$, and the sets in this union have disjoint interiors. 
In the case that $\Lb$ is a primitive substitution Delone multiset we will 
understand the term representable to mean relative to the tiles 
$T_i = (A_i,i)$, for $i\le m$, arising from the solution to the adjoint system (\ref{eq-til}). 
We call $\Lb + \Ak$ the associated tiling of $\Lb$.}
\end{defi}

\begin{defi}
{\em Let $\Lb$ be a primitive substitution Delone multiset and let
$\Pb$ be a cluster of $\Lb$. The
cluster $\Pb$ will be called {\em legal} if it is a translate of a subcluster of
$\Phi^k(x_j)$ for some $x_j \in \Lam_j$, $j \le m$ and $k \in \Z_+$.}
\end{defi}

\begin{lemma}\cite{LMS2}
Let $\Lb$ be a primitive substitution Delone multiset such that every 
$\Lb$-cluster is legal. Then $\Lb$ is repetitive. 
\end{lemma}

Not every substitution Delone multiset is representable (see \cite[Ex.\,3.12]{LMS2}), but the following theorem provides the sufficient condition for it.

\begin{theorem}\cite{LMS2}\label{legal-rep}
Let $\Lb$ be a repetitive primitive substitution Delone multiset. Then every 
$\Lb$-cluster is legal if and only if $\Lb$ is representable.
\end{theorem}

\begin{remark}
{\em In \cite[Lemma 3.2]{lawa} it is shown that if $\Lb$ is a substitution Delone multiset,
then there is a finite multiset (cluster) $\Pb \subset \Lb$ for which
$\Phi^{n-1}(\Pb) \subset \Phi^n(\Pb)$ for $n \ge 1$ and 
$\Lb = \lim_{n \to \infty} \Phi^n (\Pb)$. We call such a multiset $\Pb$ 
a {\em generating multiset}.
Note that, in order to check that every $\Lb$-cluster is legal, 
we only need to see if some cluster that contains a finite generating 
multiset for $\Lb$ is legal.}  
\end{remark}

Let $\Xi(\Tk)$ be the set of
translation vectors between $\Tk$-tiles of the same type:
\be \label{def-xi}
\Xi(\Tk) := \{x\in \R^d:\ \exists \,T,T' \in \Tk, \ T'=x+T \}.
\ee
Since $\Tk$ has the inflation symmetry with the expansive map $Q$, 
we have that
$Q\Xi(\Tk) \subset \Xi(\Tk)$.

\begin{remark} 
{\em We should be careful to distinguish between substitution Delone {\em
multisets} and substitution Delone {\em sets}. Lagarias \cite{Lag00}
considers the latter under the name of {\em self-replicating sets}.
Note that a substitution Delone set may arise from different 
substitution Delone multisets.}
\end{remark}

\subsection{Diffraction and dynamical spectra on Delone sets}
We are going to use the mathematical concept of diffraction measure developed
by Hof \cite{Hof1,Hof2}.
Given a translation-bounded measure $\nu$ on $\R^d$, let $\gamma(\nu)$ denote
its autocorrelation (assuming it is unique), that is, the vague limit
\be \label{eq-auto1}
\gamma(\nu) = \lim_{n\to \infty} \frac{1}{\Vol(F_n)} \left(
\nu|_{F_n} \ast \widetilde{\nu}|_{F_n} \right),
\ee
where $\{F_n\}_{n \ge 1}$ is a van Hove sequence
\footnote{Recall that if $f$ is a function in $\R^d$,
then $\tilde{f}$ is defined by
$\tilde{f}(x) = \overline{f(-x)}$. If $\mu$ is a measure, $\tilde{\mu}$ is
defined by $\tilde{\mu}(f) = \overline{\mu(\tilde{f})}$ for all $f \in
\Cant_0(\R^d)$.}. 
The measure $\gamma(\nu)$ is positive definite, so by Bochner's Theorem the
Fourier transform $\widehat{\gamma(\nu)}$
is a positive measure on $\R^d$, called the {\em diffraction measure} for
$\nu$.  We say  that the measure $\nu$
has {\em pure point diffraction spectrum}, if
$\widehat{\gamma(\nu)}$ is a pure point or discrete measure.
The point masses of the diffraction measure are called {\em Bragg peaks}.
For a Delone set $\Lam$ let
$$
\delta_\Lam: = \sum_{x\in \Lam} \delta_x.
$$
It is known that if $\Lam$ is a primitive 
substitution Delone set of finite local complexity, 
then $\delta_\Lam$ has a unique autocorrelation
measure $\gamma(\delta_\Lam)$(see \cite{LMS2}). We say that $\Lam$ is {\em pure point diffractive}
if the diffraction measure $\widehat{\gam(\delta_\Lam)}$ is pure discrete.

\medskip 

Let $\Lb$ be a Delone multiset and let $X_{\Lbs}$ be the collection of all 
Delone multisets each of whose
clusters is a translate of a $\Lb$-cluster. We introduce a metric
on Delone multisets in a simple variation of the standard way:
 for Delone multisets $\Lb_1$, $\Lb_2 \in X_{\Lbs}$,
\be \label{metric-multisets}
d(\Lb_1,\Lb_2) := \min\{\tilde{d}(\Lb_1,\Lb_2), 2^{-1/2}\}\, ,
\ee
where
\be 
\tilde{d}(\Lb_1,\Lb_2)
&=&\mbox{inf} \{ \e > 0 : \exists~ x,y \in B_{\e}(0), \nonumber \\ \nonumber
&  & ~~~~~~~~~~ B_{1/{\e}}(0) \cap (-x + \Lb_1) = B_{1/{\e}}(0) 
\cap (-y + \Lb_2) \}\,. 
\ee
For the proof that $d$ is a metric, see \cite{LMS1}.

Observe that $X_{\Lbs} = \overline{\{-h + \Lb : h \in \R^d \}}$ where the closure is taken in the topology induced by the metric $d$. 
The group $\R^d$ acts on $X_{\Lbs}$ by translations which are obviously
homeomorphisms, and we get a topological dynamical system $(X_{\Lbs},\R^d)$.

\begin{sloppypar}
Let $\mu$ be an ergodic invariant Borel probability
measure for the dynamical system 
$(X_{\Lbs},\R^d)$.
We consider the associated
group of unitary operators $\{U_g\}_{g\in \R^d}$ on $L^2(X_{\Lbs},\mu):$
\[U_g f(\Sk) = f(-g+\Sk).\]
A vector $\alpha =(\alpha_1,\ldots,\alpha_d) \in \R^d$ is said to be
an eigenvalue for the $\R^d$-action if there exists an eigenfunction
$f\in L^2(X_{\Lbs},\mu),$ that is, $\ f\not\equiv 0$ and
\[U_g f = e^{2 \pi i g \cdot \alpha} f,\ \ \ \mbox{for all}
\ \ g\in \R^d.\]
\end{sloppypar}

The dynamical system $(X_{\Lbs},\mu,\R^d)$ is said to have {\em pure discrete}
(or pure point) {\em spectrum} if the linear span of the
eigenfunctions is dense in $L^2(X_{\Lbs}, \mu)$.

Let $X_{\Tk} = \overline{\{-g+\Tk : g \in \R^d\}}$, where $X_{\Tk}$ carries a 
well-known topology, given analogously to (\ref{metric-multisets}) 
for $X_{\Lbs}$,
relative to which it is compact (equivalent to FLC). We have a natural action of
$\R^d$ on $X_{\Tk}$ which makes it a topological dynamical system.
The set $\{- g+\Tk:g \in \R^d \}$ is the orbit of $\Tk$.
                                                                                
Recall that a topological dynamical system is {\em uniquely ergodic}
if there is a unique invariant probability measure (which is then automatically
ergodic).
It is known (see e.g. \cite[Th.\ 2.7]{LMS1})
that for a Delone multiset $\Lb$ with FLC,
the dynamical system $(X_{\Lbs},\R^d)$ is uniquely ergodic
if and only if $\Lb$ has UCF.

\begin{theorem} \cite[Th.\ 3.2]{LMS1} \label{th-LMS1}
Suppose that a Delone multiset
$\Lb$ has FLC and UCF. Then the following are equivalent:
\begin{itemize}
\item[{\rm (i)}] $\Lb$ has pure point dynamical spectrum;
\item[{\rm (ii)}]
The measure $\nu = \sum_{i \le m} a_i \delta_{\Lambda_i}$ has pure point
diffraction spectrum, for any choice of complex numbers $(a_i)_{i\le m}$;
\item[{\rm (iii)}]
The measures $\delta_{\Lambda_i}$ have pure point diffraction spectrum,
for $i\le m$.
\end{itemize}
\end{theorem}

\section{Jordan canonical form}

Let $Q$ be a linear map from $\R^d$ to $\R^d$. We can consider $Q$ as a ($d \times d$) matrix. We discuss the matrix analysis on $Q$ that we are going to use in this paper (see \cite{HJ}).
The matrix $Q$ is similar to a matrix in the Jordan canonical form $J$, so that $Q = SJS^{-1}$ for some invertible matrix $S$ over $\C$. Suppose that $Q$ has $r$ distinct eigenvalues 
$\lambda_1, \dots, \lambda_r \in \C $. For each eigenvalue $\lambda_i, 1 \le i \le r$, there are Jordan blocks $J_{i1}(\lambda_i), \dots, J_{im_i}(\lambda_i)$ corresponding to $\lambda_i$. 
We simply write $J_{ij}$ for $J_{ij}(\lambda_i)$. 
We can decompose $J_{ij} = \lambda_i I + N $ with a matrix $\lambda_i I$ of diagonal entries and a matrix $N$ of off-diagonal entries.
For each Jordan block $J_{ij}, 1 \le j \le m_i$, we have vectors 
$e_{ij1}, \dots , e_{ijk_{ij}} \in \C^d$ such that 
\[ Qe_{ij1} = \lambda_i e_{ij1} \ \mbox{and} \ Qe_{ijl} = e_{ij(l-1)} + \lambda_i e_{ijl} 
\ \ \ \mbox{for} \ 2 \le l \le k_{ij}.
\]
For each Jordan block $J_{ij}$ and any $n \in \Z_+$, there is a simple general formula for $(J_{ij})^n$:
\[(J_{ij})^n = (\lam_i I + N)^n = \sum_{k=0}^n 
     \left( \begin{array}{c}
             n \\
             k
            \end{array} \right) \lam_i^{n-k} N^k.        
\]
We define $\left( \begin{array}{c}
             n \\
             k
            \end{array} \right) = 0$ for $n < k$.
Then for any $n \in \Z_+$,
\[ \tiny (J_{ij})^n  = \left[ \begin{array}{ccccc}
                  \lam_i^n &  \left( \begin{array}{c}
                               n \\
                               1
                              \end{array} \right) \lam_i^{n-1} &
                   \left( \begin{array}{c}
                               n \\
                               2
                              \end{array} \right) \lam_i^{n-2} &
                    \cdots  &
                   \left( \begin{array}{c}
                               n \\
                               k_{ij} - 1
                              \end{array} \right) \lam_i^{n- k_{ij} +1} \\
                   0 & \lam_i^n & \left( \begin{array}{c}
                               n \\
                               1
                              \end{array} \right) \lam_i^{n-1} & 
                     \cdots & \vdots \\
                    \vdots   & \vdots &  & & \vdots \\
                   \vdots & \vdots & \ddots  & &  \left( \begin{array}{c}
                               n \\
                               2
                              \end{array} \right) \lam_i^{n-2} \\
                    0 & 0 &  & \ddots  &  \left( \begin{array}{c}
                               n \\
                               1
                              \end{array} \right) \lam_i^{n-1} \\ \\
                    0 & 0 & \dots & \dots & \lam_i^n 
                              \end{array}  \right]
            \]

Note that $E := \{e_{ijl} \in \C^d : \ 1 \le i \le r, 1 \le j \le m_i, 1 \le l \le k_{ij}\}$ is a basis of $\C^d$.
So for any $y \in \R^d$, we can write 
\begin{equation} \label{exp1}
y = \sum_{i =1}^r \sum_{j = 1}^{m_i} \sum_{l=1}^{k_{ij}} a_{ijl}(y) e_{ijl}, 
\end{equation} 
where $a_{ijl}(y) \in \C$.

Let $\langle x, y \rangle$ be the standard inner product of $x$, $y$ in $\C^d$ and 
let $K := \mbox{max}\{k_{ij} - 1 : 1 \le i \le r, 1 \le j \le m_i \}$. 

\medskip 

\begin{lemma} \label{Polynomial-with-lambda}
Let $\alpha \in \R^d$ and $ Q : \R^d \rightarrow \R^d$ be a linear map. 
For any $n \in \Z_+$ and $w \in \R^d$ for which 
$w = \sum_{i =1}^r \sum_{j = 1}^{m_i} \sum_{l=1}^{k_{ij}} 
a_{ijl}(w) e_{ijl}$ with $a_{ijl}(w) \in \C$,  
\[ \langle \sum_{j=1}^{m_i} \sum_{l=1}^{k_{ij}} a_{ijl}(w) 
Q^n e_{ijl}, \alpha \rangle  = (P_{\alpha, w})_{i}(n) \lam_i^n \ \ \ 
\mbox{for $1 \le i \le r$}\]
and so
\[\langle Q^n w, \alpha \rangle= \sum_{i=1}^r (P_{\alpha, w})_i(n) \lam_i^n\] 
where 
$(P_{\alpha, w})_{i}$ is a polynomial over $\C$ of degree less than or equal to $K$.
\end{lemma}

\noindent
{\sc Proof.} This is standard; we provide a proof for completeness.

 We extend the linear map $Q$ from $\R^d$ to $\C^d$, i.e. 
 $Q : \C^d \rightarrow \C^d$ (just use the
same matrix). 
First note that for any $1 \le i \le r$ and $1 \le j \le m_i$, 
\begin{eqnarray*}
\lefteqn { \langle \sum_{l=1}^{k_{ij}} a_{ijl}(w) Q^n e_{ijl}, \alpha \rangle} \\
&=& \langle  a_{ij1}(w) \lam_i^n  e_{ij1}, \alpha \rangle  \\
&& + \langle a_{ij2}(w) \left( \left( \begin{array}{c}
                               n \\
                               1
                              \end{array} \right) \lam_i^{n-1} e_{ij1} + 
   \lam_i^{n} e_{ij2} \right), \alpha \rangle    \\
&&    \vdots \\
&&  + \langle a_{ijk_{ij}}(w) \left( \left( \begin{array}{c}
                               n \\
                               k_{ij} - 1
                              \end{array} \right) \lam_i^{n-k_{ij}+1} e_{ij1} + 
                   \cdots + \lam_i^{n} e_{ijk_{ij}} \right) , \alpha \rangle \,. 
\end{eqnarray*} 
Rearranging the above equation,
\begin{eqnarray*} 
\lefteqn { \langle \sum_{l=1}^{k_{ij}} a_{ijl}(w) Q^n e_{ijl}, 
\alpha \rangle} \\
&=&  \left( a_{ij1}(w) \lam_i^0 +  
      \cdots
     + a_{ijk_{ij}}(w) \left( \begin{array}{c}
                               n \\
                               k_{ij} - 1
                              \end{array} \right) \lam_i^{-k_{ij}+1} \right) 
        \langle e_{ij1}, \alpha \rangle  \lam_i^n \nonumber \\
 && + \left( a_{ij2}(w) \lam_i^0 + 
           \cdots
             + a_{ijk_{ij}}(w) \left( \begin{array}{c}
                               n \\
                               k_{ij} - 2
                              \end{array} \right) \lam_i^{-k_{ij}+2} \right) 
             \langle e_{ij2}, \alpha \rangle \lam_i^n \nonumber \\
 && \vdots  \nonumber \\
 && + \left( a_{ijk_{ij}}(w) \lam_i^0 \right) \langle e_{ijk_{ij}}, \alpha \rangle  \lam_i^n \,.
\end{eqnarray*}
Thus we get
\[ \langle \sum_{l=1}^{k_{ij}} a_{ijl}(w) Q^n e_{ijl}, \alpha \rangle  = (P_{\alpha, w})_{ij}(n) \lam_i^n, \]
where $(P_{\alpha, w})_{ij}$ is a polynomial over $\C$ of degree at most $k_{ij} - 1$. 
Then for each $1 \le i \le r$, we can write 
\be 
\langle \sum_{j=1}^{m_i} \sum_{l=1}^{k_{ij}} a_{ijl}(w) 
Q^n e_{ijl}, \alpha \rangle  = (P_{\alpha, w})_{i}(n) \lam_i^n, 
\ee
where $(P_{\alpha, w})_{i} = \sum_{j =1}^{m_i} (P_{\alpha, w})_{ij}$ is a polynomial over $\C$ of 
degree $\le K$.
Furthermore,
\[
\langle Q^n w, \alpha \rangle  = 
\langle Q^n \left( \sum_{i =1}^r \sum_{j = 1}^{m_i} \sum_{l=1}^{k_{ij}} 
a_{ijl}(w) e_{ijl} \right), \alpha \rangle  
  =   \sum_{i=1}^r (P_{\alpha, w})_i(n) \lam_i^n. 
\]
\qed

\section{Proof of the Meyer property}

The result of the following lemma is taken from \cite{Ken}.

\begin{lemma} \label{algebraic-integers}
Suppose that $L$ is a finitely generated free Abelian group in $\R^d$ such that $L$ spans $\R^d$ and $Q L \subset L$ with a linear
map $Q$. Then all eigenvalues of $Q$ are algebraic integers.
\end{lemma}

\noindent
{\sc Proof.} Let $\{{\mathrm{v}}_1, \dots, {\mathrm{v}}_n\}$ be a set of
generators for $L$.
Consider the ($d \times n$) matrix $N = [{\mathrm{v}}_1,\dots,{\mathrm{v}}_n]$.
Since $L$ spans $\R^d$, the rank of $N$ is $d$. Thus $N^T {\mathrm{x}} = {\bf 0}$ has a unique trivial solution. From the assumption of $Q L \subset L$, for each $1 \le i \le n$, we can write 
\[Q {\mathrm{v}}_i = \sum_{j=1}^n a_{ij} {\mathrm{v}}_j \ \ \ \mbox{for some} \ a_{ij} \in \Z. \]
Let $M = (a_{ij})_{n \times n}$. Then $Q N = N M^T$ and so $M N^T = N^T Q^T$. 
For any eigenvalue $\lambda$ of $Q^T$ and the corresponding eigenvector 
${\mathrm{x}}$,
\[M(N^T {\mathrm{x}}) = N^T (Q^T {\mathrm{x}}) = N^T \lambda {\mathrm{x}} = \lambda (N^T {\mathrm{x}}).\]
Since ${\mathrm{x}}$ is nonzero, $N^T {\mathrm{x}}$ is nonzero and so $\lambda$ is an eigenvalue of $M$. Since $M$ is an integer matrix, $\lambda$ is an algebraic integer. Since $Q^T$ and $Q$ have the
same eigenvalues, all eigenvalues of $Q$ are algebraic integers. \qed

\begin{cor} \label{cor-alg}
Suppose that $\Tk$ is a fixed point of a primitive substitution with expansive 
map $Q$ which has FLC. Then all eigenvalues of $Q$ are algebraic integers.
\end{cor}

\noindent
{\sc Proof.} Let $L$ be an Abelian group generated by $\Xi(\Tk)$. Since $\Tk$ has FLC, $L$ is a finitely generated free Abelian group.
From $Q \Xi(\Tk) \subset \Xi(\Tk)$ we have
$Q L \subset L$. By Lemma \ref{algebraic-integers}, 
all eigenvalues of $Q$ are algebraic integers.
\qed


\medskip

The following is a generalization of Pisot's theorem, due to K\"ornyei
\cite{Korn}.
A similar result was obtained by Mauduit \cite{Maud}. The theorem is about two equivalent conditions, 
but we state only one direction which we use later, in the special case we
need. For $x \in \R$, 
let $||x||$ denote the distance from $x$ to the nearest integer.

\begin{theorem} \cite[Th.\,1]{Korn} \label{Korn-Theorem}
Let $\lam_1, \dots, \lam_r$ be distinct algebraic numbers such that $|\lambda_i| \ge 1$, $i = 1, \dots, r$, and let
$P_1, \dots, P_r$ be nonzero polynomials with complex coefficients. 
If $\sum_{i=1}^r P_i (n) \lam_i^n$ is real for all $n$ and 
\[ \lim_{n \to \infty} || \sum_{i=1}^r P_i (n) \lam_i^n|| = 0,\]
then the following assertions are true: 
\begin{itemize} 
\item[(a)] The coefficients of $P_i$ are elements of the algebraic extension 
$\Q(\lam_i)$.
\item[(b)] If $\lam_s$ and $\lam_t$ are conjugate elements over $\Q$, and the corresponding polynomials have the form 
\[P_s (x) = \sum_{k=0}^{K_s} c_{s,k} x^k, \  P_t (x) = \sum_{k=0}^{K_t} c_{t,k} x^k, \]
then $P_s$ and $P_t$ have the same degree, $c_{w,s,k}$ and $c_{w,t,k}$ are conjugate elements over $\Q$, and for any 
isomorphism $\tau$ which is the identical mapping on 
$\Q$ and for which $\tau(\lam_s) = \lam_t$, we have
\[ \tau(c_{s,k}) = c_{t,k}, \ \ \ \mbox{for any} \ 0 \le k \le K_s = K_t.\]

\item[(c)] All the conjugates of the $\lam_i$'s not occurring in the sum 
$\sum_{i=1}^r P_i (n) \lam_i^n$ have absolute value less than one.
In other words, if $\lam'$ is a conjugate of $\lam_i$ for some $i\le r$ and
$|\lam'|\ge 1$, then $\lam' = \lam_j$ for some $j\le r$.
\end{itemize}
\end{theorem}

\begin{defi} \cite{Prag}
{\em Let $\Tk$ be a fixed point of a primitive substitution with expansive map $Q$. For each $\Tk$-tile $T$, fix a tile $\gamma T$ in the patch $\omega (T)$; 
choose $\gamma T$ with the same relative position for all tiles of the same 
type. This defines a map $\gamma : \Tk \to \Tk$ called the 
{\em tile map}. Then define the {\em control point} for a tile $T \in \Tk$ by 
\[ \{c(T)\} = \bigcap_{n=0}^{\infty} Q^{-n}(\gamma^n T).\]
}
\end{defi}

\noindent
The control points have the following properties:
\begin{itemize}
\item[(a)] $T' = T + c(T') - c(T)$, for any tiles $T, T'$ of the same type;
\item[(b)] $Q(c(T)) = c(\gamma T)$, for $T \in \Tk$.
\end{itemize}
Control points are also fixed for tiles of any tiling $\mathcal{S} \in X_{\Tk}$: they have the same relative position as in $\Tk$-tiles. 

For $n\ge 1$ let $\Tk^n := \{Q^n T:\, T\in \Tk\}$. 
By definition, if $T=(A,i)$, then $Q^nT = (Q^n A,i)$.
Thus we consider $Q^n T$ as a tile and $\Tk^n$ as a tiling. 
The tiles of $\Tk^n$ are called
{\em supertiles of level} $n$ and $\Tk^n$ is called a {\em supertiling}.
Since $\Tk$ is a fixed point of
the substitution $\om$ with expansion $Q$, we recover 
$\Tk$ by subdividing the tiles of $\Tk^n$ $n$ times. 
The control points are determined for the tiles of supertilings by
$c(Q T) = Q c(T)$. For each $T\in \Tk$ let $T^{(n)}$ be the unique
supertile of level $n$ such that $\supp(T) \subset  \supp(T^{(n)})$.

Recall that our tile-substitution $\om$ is primitive, that is,
for some $k\in \Nat$, the $k$-th power of the substitution matrix has strictly
positive entries. Then we can replace $\om$ by $\om^k$ and assume that the
substitution matrix itself is strictly positive (this does not lead to
loss of generality since a fixed point of $\om$ is also a fixed point of
$\om^k$). This means that the patch $\om(T)$ contains tiles of all types
for every $T\in \Tk$.
We can then define control points for $\Tk$-tiles choosing the tile
map $\gamma : \Tk \to \Tk$ so that for any $T \in \Tk$, the tile $\gamma T$ has the same tile
type in $\Tk$. Then for any $T, S \in \Tk$,
\[ c(\gamma T) - c(\gamma S) \in \Xi(\Tk).\]
Since $Q c(T) = c (\gamma T)$ for any $T \in \Tk$,
\be \label{control-point-translation}
Q (c(T) -  c(S)) \in \Xi(\Tk) \ \ \ \mbox{for any} \ T, S \in \Tk.
\ee

The next lemma is very close to \cite[Th.\,1.5]{Prag} and 
\cite[Lem.\,6.5]{soltil} (however, in \cite{soltil} FLC was assumed);
we provide a direct proof
for completeness.
 
\begin{lemma} \label{sum-for-translation}
Let $\Tk$ be a fixed point of a substitution with expansive map $Q$ 
and a strictly positive substitution matrix, and suppose that the control
points satisfy (\ref{control-point-translation}).
Then there exists a finite set $U$ in $\R^d$ for which $QU \subset \Xi(\Tk)$ 
and $ 0 \in U$ so that for any $T, S \in \Tk$ there exist
$N \in \N$ and $u(n), w(n)\in U,\ 0\le n \le N$, such that 
$$
c(T) - c(S) = \sum_{n=0}^N Q^n (u(n)+w(n)). 
$$
\end{lemma}

\noindent
{\sc Proof.} 
Fix any $T, S \in \Tk$ and consider the sequences of supertiles 
$T = T^{(0)} \subset T^{(1)} \subset \cdots$ and $S = S^{(0)} 
\subset S^{(1)} \subset \cdots$ defined above (to be more precise, we should
write inclusions for supports).
Fix any patch $P$ with the origin in the interior of its support.
Then there exists $N\in \Nat$ such that $T,S \in \om^N(P)$. Fix such an $N$.
Observe that
$T^{(N)} = Q^N T^{'}$ and $S^{(N)} = Q^N S^{'}$ for some 
$T^{'}, S^{'} \in P$.
We have 
\[ c(T) - c(S) = \sum_{n=0}^{N-1}\{ c(T^{(n)}) - c(T^{(n+1)})\} + 
                c(T^{(N)}) - c(S^{(N)}) -
                  \sum_{n=0}^{N-1}\{ c(S^{(n)}) - c(S^{(n+1)})\}.
\]
Note that $c(T^{(N)}) - c(S^{(N)}) = Q^N(c(T^{'}) - c(S^{'}))$ and
\begin{eqnarray*} 
c(T^{(n)}) - c(T^{(n+1)}) &=& Q^n c(T_n^{''}) - Q^n c(\gamma T_n^{'''}) \\
         &=& Q^n (c(T_n^{''}) - c(\gamma T_n^{'''})) 
\end{eqnarray*}
for some  $\Tk$-tiles $T_n^{''}, T_n^{'''}$ such that $T_n^{''} \in 
\om (T_n^{'''})$. Similarly,
\[c(S^{(n)} - c(S^{(n+1)}) = Q^n ( c(S_n^{''}) - c(\gam S_n^{'''}))\]
for  some  $\Tk$-tiles $S_n^{''}, S_n^{'''}$ such that $S_n^{''} \in      
\om (S_n^{'''})$. Thus,
\[
c(T) - c(S) = \sum_{n=0}^{N-1} Q^n\{c(T_n^{''}) - c(\gamma T_n^{'''}) -
(c(S_n^{''}) - c(\gam S_n^{'''}))\} + Q^N (c(T^{'}) - c(S^{'})).
\]
Observe that there are finitely many possibilities for 
$c(T_n^{''}) - c(\gamma T_n^{'''})$, $c(S_n^{''}) - c(\gam S_n^{'''})$, and
$c(T^{'}) - c(S^{'})$ (for the first two differences
it suffices to consider all the cases for which $T_n^{'''}$ and $S_n^{'''}$ are prototiles, $T_n^{''} \in \om (T_n^{'''})$ and $S_n^{''} \in      
\om (S_n^{'''})$ ).
Thus, we obtained the desired representation, in view of 
(\ref{control-point-translation}).
\qed

\begin{theorem} \cite[Th.\,4.3]{soltil} \label{translation-vector-formula}
Let $\Tk$ be a repetitive fixed point of a primitive substitution with expansive map $Q$ which has FLC.
If $\alpha \in \R^d$ is an eigenvalue for $(X_{\Tk}, \R^d, \mu)$, 
then for any $x \in \Xi(\Tk)$ we have 
$||\langle Q^n x, \alpha \rangle|| \stackrel{n \to \infty}{\rightarrow} 0$.
\end{theorem}

In \cite{soltil} it was assumed that the expansive map $Q$ is diagonalizable
over $\C$, but the proof works in full generality.

Let $\mathcal{M} = \{(c(T) - c(S)) - (c(T') - c(S')) \,:
\ T, S, T', S' \in \Tk \} \subset \R^d$.
Combining Lemma \ref{sum-for-translation} and Th.\,\ref{translation-vector-formula}, we obtain the following corollary.

\begin{cor} \label{sum-of-finite-elements}
Let $\Tk$ be a fixed point of a substitution 
with expansive map $Q$ and a strictly positive substitution matrix
which has FLC. Suppose that the control points satisfy
(\ref{control-point-translation}).
Let $\alpha \in \R^d$ be an eigenvalue for $(X_{\Tk}, \R^d, \mu)$. 
Then there exists a finite subset 
$W$ in $\R^d$ independent to the choice of $\alpha$ for which 
$$
||\langle Q^n w, \alpha \rangle|| \stackrel{n \to \infty}{\rightarrow} 0 
\ \ \ \mbox{for any} \ w \in W,
$$
and for any $y \in \mathcal{M}$, 
there exist $N\in \Nat$ and $w(n) \in W$,
$0 \le n\le N$, such that
$$
y = \sum_{n=0}^N Q^n w(n).
$$
\end{cor}

\begin{prop} \label{tiling-meyer}
Let $\Tk$ be a fixed point of a substitution
with expansive map $Q$ and a strictly positive substitution matrix
which has FLC. Suppose that the control points satisfy
(\ref{control-point-translation}) and
the set of eigenvalues for $(X_{\Tk}, \R^d, \mu)$ 
is relatively dense.
Then $\{c(T) - c(S) \,: \ T, S \in \Tk \}$ is uniformly discrete, that is,
$\{c(T) \,: \ T \in \Tk \}$ is a Meyer set.
\end{prop}

\noindent
{\sc Proof.} Since the set of eigenvalues is relatively 
dense, there exist 
eigenvalues $\alpha_1, \dots, \alpha_d$ for $(X_{\Tk}, \R^d, \mu)$ such that 
for any $0 \neq y \in \R^d$,
\[\langle y, \alpha_t \rangle \neq 0 \ \ \ \mbox{for some} \ 1 \le t \le d.\]
We define a norm $|||\cdot|||$ on $\R^d$ in terms of the expansion
(\ref{exp1}):
\[|||y||| = 
  ||| \sum_{i =1}^r \sum_{j = 1}^{m_i} \sum_{l=1}^{k_{ij}} a_{ijl}(y) e_{ijl} ||| 
 = \sum_{t=1}^d \left( \sum_{i =1}^r \sum_{j = 1}^{m_i} \sum_{l=1}^{k_{ij}} 
             |a_{ijl}(y) \langle e_{ijl}, \alpha_t \rangle | \right). \]
For any $\alpha \in \{\alpha_1,\ldots,\alpha_d\}$ we have
\begin{equation}\label{eka1}
\langle y,\alpha \rangle = \sum_{i=1}^r T_{y,\alpha,i},\ \ \ \mbox{where}\ \ 
T_{y,\alpha,i} = \sum_{j = 1}^{m_i} \sum_{l=1}^{k_{ij}} a_{ijl}(y)
\langle e_{ijl}, \alpha \rangle.
\end{equation}
Clearly,
\begin{equation} \label{esq1}
|||y||| \ge \sum_{i=1}^r |T_{y,\alpha,i}|.
\end{equation}

From Cor.\,\ref{sum-of-finite-elements} we know that any $y \in 
\mathcal{M}$
can be represented in terms of elements of $W$ so that 
$y = \sum_{n=0}^N Q^n w(n)$ for some positive integer $N$, where $w(n) \in W$ 
for any $0 \le n \le N$. 
Let $w_1, \dots, w_R$ be all the elements of $W$. 
We can rearrange the sum to write
\begin{equation} \label{eka2}
y = \sum_{p=1}^R \sum_{n \in \Nk_p} Q^{n} w_p, 
\end{equation}
where $\{\Nk_1, \dots, \Nk_R\}$ is 
a partition of $\{0, 1, \dots, N\}$ such that 
$w(n) = w_p$ if and only if $n \in \Nk_p$ for $0 \le n \le N$.
For any $\alpha \in \{\alpha_1, \dots, \alpha_d\}$ and 
$w = \sum_{i =1}^r \sum_{j = 1}^{m_i} \sum_{l=1}^{k_{ij}} 
a_{ijl}(w) e_{ijl} \in W$, 
by Lemma \ref{Polynomial-with-lambda} we have
\be \label{new4}
\langle \sum_{j=1}^{m_i} \sum_{l=1}^{k_{ij}} a_{ijl}(w) 
Q^n e_{ijl}, \alpha \rangle  = (P_{\alpha, w})_{i}(n) \lam_i^n \ \ \ 
\mbox{for $1 \le i \le r$} 
\ee
and
\begin{equation} \label{new2}
\langle Q^n w, \alpha \rangle 
= \sum_{i=1}^r (P_{\alpha,w})_i (n) \lam_i^n,  
\end{equation}
where $(P_{\alpha,w})_i$ is a polynomial over $\C$
of degree less than or equal to $K$.
Comparing (\ref{eka1}) and (\ref{eka2}) and noting that $Q^n e_{ijl}$ is in 
the subspace of $\C^d$ spanned by $e_{ij1}, \dots, e_{ijk_{ij}}$, we obtain 
\[
T_{y,\alpha,i} = \sum_{p=1}^R \sum_{n \in \Nk_p} \langle \sum_{j=1}^{m_i} \sum_{l=1}^{k_{ij}} 
a_{ijl}(w_p) Q^n e_{ijl}, \alpha \rangle \ \ \ \mbox{for $1 \le i \le r$}.
\]
From (\ref{new4}), we get 
\begin{equation} \label{eka3}
T_{y,\alpha,i} = \sum_{p=1}^R \sum_{n \in \Nk_p}
(P_{\alpha,w_p})_i (n) \lam_i^n \ \ \ \mbox{for $1 \le i \le r$}.
\end{equation}
Note that $||\langle Q^n w, \alpha \rangle|| \stackrel{n \to \infty}{\to} 0$ 
by Cor.\,\ref{sum-of-finite-elements}, and for any $1 \le i \le r$, $\lam_i$ is 
an algebraic integer by 
Cor.\,\ref{cor-alg}, with $|\lam_i| > 1$ by expansiveness of $Q$.
Therefore, by Th.\,\ref{Korn-Theorem}, for any $w \in W$ and $1 \le i \le r$ we have
\be \label{Pi-polynomial}
(P_{\alpha,w})_i (n) = \sum_{k=0}^K  (c_{\alpha,w,i,k})n^k,
\ee 
where
$c_{\alpha,w,i,k} \in \Q(\lambda_i)$, 
and every conjugate $\lam$ of $\lam_i$, with $|\lam|\ge 1$, occurs in the
right-hand side of (\ref{new2}), that is, $\lam=\lam_j$ for some $j\le r$.
Moreover, in this case
\[ c_{\alpha,w,j,k} = \tau_{ij}(c_{\alpha,w,i,k}) \ \ \ \mbox{for any 
$0 \le k \le K$,}\]
where $\tau_{ij}:\ \Q(\lam_i) \to \Q(\lam_j)$
is an isomorphism which is identical on $\Q$ such that
$\tau_{ij}(\lambda_i) = \lambda_j$.
Since all $\lambda_i$ are algebraic integers, we have
$$
\Q(\lambda_i) = \Q[\lam_i] =
\{a_0 + a_1 \lambda_i + \cdots + a_{s_i-1} \lambda_i^{s_i-1} :\ 
a_n \in \Q, \ 0 \le n \le s_i-1\},
$$
where $s_i$ is the degree of the minimal 
polynomial of $\lambda_i$ over $\Q$.
There are finitely many numbers $c_{\alpha,w,i,k}$, so we can find
a positive integer $b$ such that
$$
b c_{\alpha,w,i,k} \in \Z[\lam_i],\ \ \forall\, \alpha\in 
\{\alpha_1,\ldots,\alpha_d\},\ \forall\,w \in W, \ \forall\,i\le r,\
\forall\, k\le K.
$$
That is, there exist polynomials $g_{\alpha,w,i,k}(x)$ with integer coefficients
such that
\begin{equation} \label{new3}
b c_{\alpha,w,i,k} = g_{\alpha,w,i,k}(\lam_i)
\end{equation}
and
$$
\lam_i,\,\lam_j\ \mbox{are conjugates} \ \Rightarrow\ g_{\alpha,w,i,k}(x)
= g_{\alpha,w,j,k}(x).
$$
Let
\begin{equation} \label{defC1}
C_1 := \max\{|g_{\alpha,w,i,k}(x)|:\ |x|\le 1,\,\alpha\in
\{\alpha_1,\ldots,\alpha_d\},\,w\in W,\,i\le r,\,
k \le K\}.
\end{equation}
Note that $C_1 < \infty$.

Now fix $0 \ne y \in \mathcal{M}$ and choose $\alpha \in \{\alpha_1, \ldots,
\alpha_d\}$ such that $\langle y, \alpha \rangle \neq 0$. Then fix 
$1 \le i \le r$ such that $T_{y,\alpha,i} \neq 0$, see (\ref{eka1}).
Consider a polynomial $S(x) = S_{y,\alpha,i}(x) \in \Z[x]$ given by
\be \label{S-polynomial}
 S(x) = \sum_{p=1}^R \sum_{n \in \Nk_p} \sum_{k=0}^K g_{\alpha,w,i,k}(x)
n^k x^n.
\ee 
In view of (\ref{eka3}), (\ref{Pi-polynomial}), (\ref{new3}),
and (\ref{S-polynomial}),
\begin{equation} \label{new1}
S(\lam_i) = b T_{y,\alpha,i}.
\end{equation}
Let $\mathcal{H}_i = 
\{\mbox{all conjugates $\lam$ of} \ \lam_i : |\lam| \ge 1 \}$ and 
$\mathcal{G}_i = \{\mbox{all conjugates $\lam$ of} \ \lam_i\}$.
By Th.\,\ref{Korn-Theorem}(c) we have $\Hk_i \subset \{\lam_1,\ldots,\lam_r\}$
and
\[\lam_j\in \Hk_i\ \Rightarrow\ S(\lam_j) = \tau_{ij}(S(\lam_i)).
\]
On the other hand,
for any $\lambda \in \mathcal{G}_i \backslash \mathcal{H}_i$,
\begin{eqnarray*} 
|S(\lam)| 
\le  C_1 \sum_{k=0}^K \sum_{n=0}^{\infty} |n^k \lam^n|
\end{eqnarray*}
where $C_1$ was defined in (\ref{defC1}).
Since $\sum_{n=0}^{\infty} n^k \lam^n$ converges absolutely for any 
$|\lam| < 1$ and $0 \le k \le K$, there exists a constant $C_2 >0$,
independent of $y$, $\alpha$, $i$, such that 
$$
|S(\lam)| < C_2  \ \ \ \mbox{for any}\ \lam \in \Gk_i\setminus \Hk_i.
$$

Now observe that
$$
\Phi:=\prod_{\lam \in \mathcal{G}_i} S(\lam) \in \Z,
$$
since $S$ is a polynomial over $\Z$ and the product
is symmetric under permutations of the conjugates of $\lam_i$.
On the other hand, $\Phi\ne 0$, since
$S(\lambda_i) = b T_{y,\alpha,i} \ne 0$ and therefore,
$S(\lambda)= \tau(S(\lambda_i)) \ne 0$ where $\tau:\, \Q(\lam_i)\to \Q(\lam)$
is an  
isomorphism satisfying $\tau(\lam_i)= \lam$, for $\lam\in \Gk_i$.
Therefore, $|\Phi|\ge 1$, hence
\begin{eqnarray}
 \prod_{\lam \in \mathcal{H}_i} |S(\lam)| 
\ge  \frac{1}{\prod_{\lam \in \mathcal{G}_i \backslash \mathcal{H}_i} 
|S(\lam)|}\,. 
\end{eqnarray}
Note that 
$$\prod_{\lam \in \mathcal{G}_i \backslash \mathcal{H}_i} 
|S(\lam)| \le (C_2)^L \ \ \  
\mbox{where} \ \ L = \# (\mathcal{G}_i \backslash \mathcal{H}_i). $$
Let $H = \#{\mathcal{H}_i}$. We obtain
\[
\left( \sum_{\lam \in \mathcal{H}_i} |S(\lam)| \right)^{H}
\ge 
\prod_{\lam \in \mathcal{H}_i} |S(\lam)| \ge (C_2)^{-L},
\]
and, in view of (\ref{esq1}) and (\ref{new1}),
\be \label{final-inequality}
 |||y||| \ge \sum_{\lam \in \mathcal{H}_i} |T_{y,\alpha,i}| =
\frac{1}{b}\sum_{\lam \in \mathcal{H}_i} |S(\lam)|
 \ge \frac{1}{b} (C_2)^{-L/H}.  
 \ee
Thus,
$\{|||y||| : y \in \mathcal{M}, y \neq 0\}$ 
has a uniform positive lower bound.
Since all norms in $\R^d$ are equivalent, the set
$\{c(T) - c(S) \,: \ T, S \in \Tk \}$  is uniformly discrete in 
the Euclidean norm. This completes the proof of the proposition.
\qed

\medskip

\begin{cor} \label{relativelyDense-Meyer}
Let $\Lb$ be a primitive substitution Delone multiset with expansion $Q$ for which every $\Lb$-cluster is legal and $\Lb$ has FLC.
If the set of eigenvalues for $(X_{\Lbs}, \R^d, \mu)$ is relatively dense, then
$\Lambda = \cup_{i \le m} \Lambda_i$ is a Meyer set.
\end{cor}

\noindent
{\sc Proof.} Since $\Lb$ is representable by Th.\,\ref{legal-rep}, 
we have that $\Tk:=\Lb + \Ak$ is a repetitive tiling which has FLC and
is a fixed point of a
primitive substitution $\om$ with expansion $Q$.
Since $(X_{\Lbs}, \R^d, \mu)$ and $(X_{\Tk}, \R^d, \mu)$ are 
topologically conjugate (see \cite[Lem.\,3.10]{LMS2}), 
the set of eigenvalues for $(X_\Tk, 
\R^d, \mu)$ is relatively dense. 
The substitution $\om$ is primitive, so we can find $k\in \Nat$ such that 
$\om^k$ has a strictly positive substitution matrix. Then we can consider
$\Tk$ as a fixed point of $\om^k$ with expansive map $Q^k$.
We can choose control
points for $\Tk$ to satisfy (\ref{control-point-translation}), with
$Q$ replaced by $Q^k$. Then Proposition~\ref{tiling-meyer} applies, and
we obtain that $\mathcal{L} - \mathcal{L}$ is uniformly discrete, 
where $\mathcal{L} := \{c(T) : T \in \Lb + \Ak \}$.

Then for each $i \le m$, $\Lambda_i \subset a_i + \mathcal{L}$ for
some $a_i \in \R^d$ and $\Lambda = \cup_{i \le m} \Lambda_i \subset
F + \mathcal{L}$ for some finite set $F$ of $\R^d$. So $\Lambda - \Lambda
\subset (F - F) + \mathcal{L} - \mathcal{L}$. Since
$(F - F) + \mathcal{L} - \mathcal{L}$ is uniformly discrete, 
$\Lambda$ is a Meyer set. 
\qed

\medskip

\begin{lemma} \label{purePoint-eigen-span}
Let $\Lb$ be a Delone multiset in $\R^d$.
Suppose that $(X_{\Lbs}, \R^d, \mu)$ has a pure point dynamical spectrum.
Then the eigenvalues for the dynamical system $(X_{\Lbs}, \R^d, \mu)$ span $\R^d$.
\end{lemma}

\noindent
{\sc Proof.} Suppose that there is a non-zero $x \in \R^d$ such that $\langle x, \alpha \rangle = 0$ for any 
eigenvalue $\alpha$ for $(X_{\Lbs}, \R^d, \mu)$. We take $x \in \R^d$ with small norm so that 
$a + x \notin \Lambda$ for all $a \in \Lambda = \cup_{i \le m} \Lambda_i$. For an eigenfunction $f_{\alpha}$ corresponding to the eigenvalue 
$\alpha$,
\[
f_{\alpha}(\Lb' - x) = e^{2 \pi i \langle x, \alpha \rangle}f_{\alpha}(\Lb') = f_{\alpha}(\Lb'), 
\ \ \ \mbox{for $\mu$-a.e.} \ \Lb' \in X_{\Lbs}. \]
For any $f \in L^2(X_{\Lbs}, \mu)$, $f = \sum_{n = 1}^{\infty} f_{\alpha_n}$, where $f_{\alpha_n}$'s are eigenfunctions.
We denote the norm in $L^2(X_{\Lbs}, \mu)$ by $\|\cdot\|_2$.
For any $\epsilon >0$, there is $N \in \N$ such that 
\begin{eqnarray*}
|| f(\cdot -x) - f||_2 &\le& ||f(\cdot - x) - \sum_{n=1}^N f_{\alpha_n}(\cdot - x)||_2 + 
                          ||\sum_{n=1}^N f_{\alpha_n}(\cdot - x) - f||_2 \\
            &\le& ||f(\cdot - x) - \sum_{n=1}^N f_{\alpha_n}(\cdot - x)||_2 + 
                          ||\sum_{n=1}^N f_{\alpha_n} - f||_2 \\
          &\le& 2 \epsilon.
\end{eqnarray*}
So $f(\Lb' -x) = f(\Lb')$ for $\mu$-a.e. $\Lb' \in X_{\Lbs}$. 
Note that $\Lb\ne \Lb-x$ by the choice of
$x$. Therefore, we can
choose $\eps > 0$ such that the $\eps$-neighborhood of $\Lb$
 and its translation by $x$
are disjoint, by the continuity of the action. Consider $f$ to be the
characteristic function of the $\eps$-neighborhood of $\Lb$. We have 
$f(\Lb')
=1$ but $f(\Lb'-x)$ = 0 for all $\Lb'$ in this neighborhood, 
which is a contradiction.
\qed

\medskip

Noticing that every integral linear combination of the eigenvalues for 
$(X_{\Lbs}, \R^d, \mu)$ is also an eigenvalue for the dynamical system, 
from Cor.\,\ref{relativelyDense-Meyer} and Lemma \ref{purePoint-eigen-span} we 
get the following theorem.

\begin{theorem} \label{Ppd-onPointset-Meyer}
Let $\Lb$ be a primitive substitution Delone multiset with expansion $Q$ for which every $\Lb$-cluster is legal and $\Lb$ has FLC.
Suppose that $(X_{\Lbs}, \R^d, \mu)$ has a pure point dynamical spectrum. Then 
$\Lambda = \cup_{i \le m} \Lambda_i$ is a Meyer set.
\end{theorem}

\begin{theorem} \cite{Nik} \label{Nik-theorem}
If $\Lambda$ is a Meyer set and its autocorrelation exists with respect to a van Hove sequence, then the set of Bragg peaks is relatively dense.
\end{theorem}

\begin{lemma} \label{Nik-lem}
Let $\Lb$ be a Delone multiset for which $\Lb$ has FLC and UCF.
If the union of the Bragg peaks of the sets $\Lambda_j$, $1 \le j \le m$, 
is relatively dense, then the set of eigenvalues for $(X_{\Lbs}, \R^d, \mu)$ 
is relatively dense.
\end{lemma}

\noindent
{\sc Proof.} This follows from \cite[Lem.\,3.4]{LMS1},
which was essentially taken from \cite{Dw,Hof2}. 
We refer \cite{LMS1} for more details.

It is enough to show that every Bragg peak of any set 
$\Lambda_j$ is an eigenvalue for $(X_{\Lbs}, \R^d, \mu)$.  
Let $\gamma = \gamma(\delta_{\Lam_j})$ denote the autocorrelation of 
$\delta_{\Lam_j}$ given by (\ref{eq-auto1}).
Let $\om \in \Cant_0(\R^d)$,
that is, $\om$ is continuous and has compact support.
We define
$$
f_{j,\om}(\Lb') := (\om\ast\delta_{{\Lam_j'}})(0)\ \ \ 
\mbox{for}\ \Lb' = (\Lambda_i')_{i\le m} \in X_{\Lbs}.
$$
Denote by $\gam_{\om,{\Lamt_j}}$ the autocorrelation of
$\om\ast\delta_{{\Lamt_j}}$.
Then 
$
\gam_{\om,{\Lamt_j}} = (\om \ast \widetilde{\om}) \ast \gamma
$
and, therefore,
$
\widehat{\gamma_{\om,{\Lamt_j}}} = |\widehat{\om}|^2 \widehat{\gamma}.$
By Lemma 3.4 in \cite{LMS1} we note that $$
\sigma_{f_{j,\om}} = \widehat{\gamma_{\om,{\Lamt_j}}},$$
where $\sigma_{f_{j,\om}}$ is the spectral measure corresponding to $f_{j,\om}$.
(In \cite[Lem.\,3.4]{LMS1} we considered the measure $\nu = \sum_{i\le m}
a_i \delta_{\Lam_i}$; here we take $a_i = \delta_{ij}$.)
If $\alpha$ is a Bragg peak of $\Lam_j$, then $\widehat{\gamma}(\alpha) >0$.
We can certainly find $\om\in \Cant_0(\R^d)$ such that $\widehat{\om}(\alpha)
\ne 0$, and  then $\sigma_{f_{j,\om}}(\alpha)>0$. Thus, the spectral measure
corresponding to some $L^2$ function has a point mass at $\alpha$, and this
implies that $\alpha$ is an eigenvalue for
the group of unitary operators (see e.g.\
\cite{Weidmann}); we conclude that
$\alpha$ is an eigenvalue for $(X_{\Lbs}, \R^d, \mu)$.
\qed

\medskip

Combining the results above we obtain the following equivalences.

\begin{theorem} \label{Bragg-relativelyDense-Meyer}
Let $\Lb$ be a primitive substitution Delone multiset with expansion $Q$ for which every $\Lb$-cluster is legal and $\Lb$ has FLC. Then the following are equivalent;
\begin{itemize}
\item[(i)] The set of Bragg peaks for each $\Lambda_j$ is relatively dense. 
\item[(ii)] The union of Bragg peaks of $\Lam_j$, $1\le j \le m$, 
is relatively dense.
\item[(iii)] 
The set of eigenvalues for $(X_{\Lbs}, \R^d, \mu)$ is relatively dense.
\item[(iv)] $\Lambda = \cup_{j \le m} \Lambda_j$ is a Meyer set.
\end{itemize}
\end{theorem}

{\em Proof.} (i) $\Rightarrow$ (ii) is trivial; (ii) $\Rightarrow$ (iii) is
Lem.~\ref{Nik-lem}, (iii) $\Rightarrow$ (iv) is 
Cor.~\ref{relativelyDense-Meyer}. Finally, (iv) $\Rightarrow$ (i) follows
by Strungaru's Th.\,\ref{Nik-theorem}. 
Note that each $\Lam_j$ is a
Meyer set, since $\Lam_j - \Lam_j$ is uniformly discrete and $\Lam_j$ is
a Delone set.  We apply Th.\,\ref{Nik-theorem} to each $\Lam_j$.
(It is known that
a primitive substitution Delone multiset for which every $\Lb$-cluster is 
legal has UCF, see e.g. \cite{LMS2}, hence for every $\Lam_j$ there exists
unique autocorrelation.)
\qed

\medskip

This theorem readily shows Theorem \ref{th-main} and Corollary \ref{cor1} in the introduction.

\medskip

{\bf Acknowledgment.} We are grateful to the referees for many helpful
comments.

\end{document}